\documentclass{article}

\usepackage{latexsym}
\usepackage{a4wide}
\usepackage{amscd}
\usepackage{graphics}
\usepackage{amsmath}
\usepackage{amssymb}
\input xy
\xyoption{all}

\newcommand{\N}{\mathbb{N}}                     
\newcommand{\Z}{\mathbb{Z}}                     
\newcommand{\R}{\mathbb{R}}                     
\newcommand{\set}[2]{\left\{{#1}\mid{#2}\right\}}       
\newcommand{\proof}{{\sl Proof.}\hspace{5pt}}   
\newcommand{\qed}{\hfill $\Box$ \bigskip}       
\newcommand{\ind}{\mathrm{ind\,}}               
\newcommand{\codim}{\mathrm{codim}}           
\newcommand{\supp}{\mathrm{supp\,}}             

\newcommand{\crit}{\mathrm{crit}}

\newcommand{\grad}{\mathrm{grad\,}}

\newtheorem{thm}{\sc Theorem}[section]      

\newtheorem{lem}[thm]{\sc Lemma}            
\newtheorem{prop}[thm]{\sc Proposition}     

\title{When the Morse index is infinite}
\author{Alberto Abbondandolo and Pietro Majer}
\date{March 31, 2004}

\begin{document}

\maketitle

\begin{abstract}
Let $f$ be a smooth Morse function on an infinite dimensional
separable Hilbert manifold, all of whose critical points have infinite
Morse index and co-index. For any critical point $x$ choose an integer
$a(x)$ arbitrarily. Then there exists a Riemannian structure on $M$
such that the corresponding gradient flow of $f$ has the following
property: for any pair of critical points $x,y$, the unstable manifold
of $x$ and the stable manifold of $y$ have a transverse intersection
of dimension $a(x)-a(y)$.
\end{abstract}

\renewcommand{\theenumi}{\roman{enumi}}
\renewcommand{\labelenumi}{(\theenumi)}

\section*{Introduction}

Let $f$ be a smooth Morse function on an infinite dimensional
separable Hilbert
manifold\footnote{By this we mean a paracompact separable space with a smooth
atlas of charts taking values in an infinite dimensional separable
real Hilbert space.} $M$. Let us denote by $\crit(f)$ the set of its
critical points, and let us assume that each $x\in \crit(f)$ has
finite Morse index $i(x)$. A Riemannian structure $g$ on $M$
determines the vector field $-\grad f$, whose local flow $\phi_t: M
\rightarrow M$ has the critical points of $f$ as rest points. The
Morse condition is translated into the fact that these rest points are
hyperbolic. Their unstable and stable manifolds,
\begin{eqnarray*}
W^u(x;f,g) = W^u(x) = \set{p\in M}{\lim_{t\rightarrow -\infty}
  \phi_t(p) = x}, \\  W^s(x;f,g) = W^s(x) = \set{p\in
  M}{\lim_{t\rightarrow +\infty} \phi_t(p) = x},
\end{eqnarray*}
are submanifolds with $\dim W^u(x) = \codim\, W^s(x) =
i(x)$. Moreover for a generic choice of the Riemannian structure $g$, 
$W^u(x)$ and $W^s(y)$ intersect
transversally for every pair of critical points $x,y$.
It follows that $W^u(x)\cap W^s(y)$ - if non-empty - is a submanifold 
of dimension
\begin{equation}
\label{launo}
\dim W^u(x) \cap W^s(y) = i(x) - i(y).
\end{equation}

If the critical points have infinite Morse index and co-index, their
unstable and stable manifolds are infinite dimensional and so could be
their intersections. However, there are situations in which these
intersections are indeed finite dimensional, and one can
associate an integer $i_{\mathrm{rel}} (x)$ to each critical point $x$
in such a way that
\begin{equation}
\label{ladue}
\dim W^u(x) \cap W^s(y) = i_{\mathrm{rel}}(x) - i_{\mathrm{rel}}(y).
\end{equation}
This fact allows to develop the analogue of Morse theory for such
functions: it is the case of Floer homology\footnote{Actually, in
Floer homology one does not consider the gradient flow with respect to
a Riemannian metric on a Hilbert manifold, but one uses the gradient
equation with respect to an inner product which is not complete on the
space where the function is smooth. Therefore, one does not obtain a
local flow and there are no stable and unstable manifolds. 
However, the space of solutions connecting two critical points is a
finite dimensional manifold.}
for the Hamiltonian action
functional on the space of loops on certain symplectic manifolds
\cite{flo88a, flo88d, flo89a}, of functionals on Hilbert spaces
somehow compatible with either a fixed splitting into two closed
linear subspaces \cite{ama01}, or with a closed linear subspace and a 
fixed flag of finite dimensional linear subspaces
\cite{szu92,ks97,gip99}, and more generally of functionals on Hilbert
manifolds which are compatible with a fixed subbundle of the tangent
bundle \cite{ama03c} or with a polarization \cite{cjs95}.

The aim of this note is to prove the following result:

\medskip

\noindent {\sc Theorem.}
{\em Let $f$ be a smooth Morse function on the Hilbert manifold $M$, all of
whose critical points have infinite Morse index and co-index. Let}
\[
a : \crit(f) \rightarrow \Z
\]
{\em be an {\em arbitrary} function. 
Then there exists a Riemannian structure $g$
on $M$ such that for every $x,y\in \crit(f)$ the intersection
$W^u(x;f,g)\cap W^s(y;f,g)$ is transverse and - if it is non-empty -
it is a submanifold of dimension $a(x) - a(y)$.}

\noindent {\em 
Moreover, if $\set{(x_i,y_i)}{i=1,\dots,n}$ is a finite set of pairs
of critical points such that $a(x_i)>a(y_i)$ and there is a smooth curve
$u_i:[0,1]\rightarrow M$, $u_i(0)=x_i$, $u_i(1)=y_i$, $i=1,\dots, n$,
such that $Df(u_i)[u_i^{\prime}]<0$ on $]0,1[$, the Riemannian
structure $g$ can be chosen in such a way that $W^u(x_i;f,g)\cap
W^s(y_i;f,g)\neq \emptyset$ for every $i=1,\dots,n$.}

\noindent{\em
Finally, if $g_0$ is any Riemannian structure on $M$, the Riemannian
structure $g$ can be chosen to be equivalent to $g_0$, meaning that
there is $c>0$ for which}
\[
\frac{1}{c} \, g_0(\xi,\xi) \leq g(\xi,\xi) \leq c \, g_0(\xi,\xi) \quad
\forall \xi \in TM.
\]

\medskip

Therefore the situation is completely different from the case of
finite Morse indices, where the identity (\ref{launo}) holds for every
choice of the Riemannian structure and Morse theory only depends on
the pair $(M,f)$. In the case of infinite Morse indices and co-indices
no Morse theory based just on the pair $(M,f)$ can possibly exist, and
extra structures, such as the ones used in the above mentioned papers,
are really needed. 

Notice that if the function $a$ is chosen to be constant, the above
theorem produces a Riemannian structure $g$ for which the unstable and
stable manifolds of two distinct critical points never meet. The
statement about the possibility of requiring some intersections to be
non-empty says that choosing a non-constant $a$ one can actually
produce non-trivial intersections.

Finally, the possibility of choosing $g$ to be equivalent to a given
Riemannian structure says that from the point of view of completeness
and compactness the Riemannian structure $g$ is not worse than a
preferred Riemannian structure $g_0$ we might dispose of: if $(M,g_0)$
is complete so is $(M,g)$, if $f$ satisfies the Palais-Smale
condition with respect to $g_0$ it will also satisfy it with respect
to $g$.

The idea of the proof is the following. After choosing a suitable Riemannian
structure in a neighborhood of the critical points, we can find an
integrable subbundle $\mathcal{V}$ of $TM$ having infinite dimension
and codimension and such that for every critical point $x$ the
negative eigenspace of the Hessian of $f$ at $x$ is a {\em compact
perturbation} of $\mathcal{V}(x)$, of {\em relative dimension} $a(x)$. The
existence of $\mathcal{V}$ follows from a well known result of Eells and
Elworthy \cite{ee70} stating in particular
that every separable infinite dimensional Hilbert
manifold can be smoothly embedded as an open subset of a Hilbert space
(hence its tangent bundle has many integrable subbundles) and from
Kuiper theorem \cite{kui65} stating that the general linear group of
an infinite dimensional Hilbert space is contractible.
The fact that $\mathcal{V}$ is integrable allows us to construct a vector
field $X$ on $M$ having $f$ as a Lyapunov function, and whose local
flow {\em essentially preserves} $\mathcal{V}$. As it is shown in
\cite{ama03c}, these facts imply that the unstable and stable
manifolds of two rest points $x,y$ of the vector field $X$ have a
{\em Fredholm intersection} of index $a(x)-a(y)$. The conclusion follows
from the fact that $X$ is the negative gradient of $f$ with respect to
a suitable Riemannian structure on $M$, which can be perturbed in
order to have transverse intersections.

\section{Flows on a Hilbert space which essentially preserve a closed
  linear subspace}  

Let $H$ be the infinite dimensional separable real Hilbert space, with 
inner product $\langle \cdot, \cdot \rangle$ and norm $|\cdot|$. The
operator norm will be denoted by $\|\cdot\|$. The set of all closed
linear subspaces of $H$ will be denoted by $Gr(H)$, and
$Gr_{\infty,\infty}(H)$ will denote the subset of all those subspaces
having infinite dimension and codimension. 

We start by recalling some basic definitions and some results from
\cite{ama03c}. A pair $(V,W)$ of closed linear subspaces of $H$ is
said a {\em Fredholm pair} if $\dim V\cap W<\infty$ and
$\codim (V+W)<\infty$, in which case the number $\ind (V,W) = \dim
V\cap W - \codim (V+W)$ is said the {\em Fredholm index} of
$(V,W)$. Two submanifolds of $H$ have a {\em Fredholm
intersection} if at every point of their intersection the pair
consisting of their tangent spaces is Fredholm. 

Given $V,W$ two closed linear subspaces of the Hilbert space $H$, $V$
is said to be a {\em compact perturbation} of $W$ if the difference of
the orthogonal projections onto $V$ and onto $W$ is a compact
operator, or equivalently if the operators
\begin{equation}
\label{ecom}
P_V P_{W^{\perp}} \quad \mbox{and} \quad P_{W} P_{V^{\perp}}
\end{equation}
are compact. In this case, $(V,W^{\perp})$ is a Fredholm pair and its
index is called the {\em relative dimension} of $V$ with respect to
$W$, and it is denoted by $\dim (V,W)$.

Let $M$ be an open subset of $H$, and
let $X$ be a smooth vector field on $M$, having
only hyperbolic rest points (this means that the differential of $X$ at
every rest point does not have purely imaginary spectrum). Let
$W$ be a closed linear subspace of $H$, with orthogonal projector $P=P_W$, 
and consider the following compatibility conditions between $X$ and $W$:

\begin{description}

\item[(C1)] For every rest point $x$ of $X$, the positive eigenspace
  of the differential of $X$ at $x$ is a compact perturbation of
  $W$.

\item[(C2)] For every $u\in H$, the operator $[DX(u),P]P = (I-P) DX(u)
  P$ is compact. 

\end{description}

By (C1), the relative dimension of the positive eigenspace of $DX(x)$ with
respect to $W$ is a well defined integer, which is 
said the {\em relative Morse index} of $x$, denoted by $i(x,W)$.
If $\phi_t$ denotes the local flow determined by the vector field $X$,
condition (C2) is equivalent to the fact that $\phi_t$ {\em essentially
preserves} $W$, in the sense that $D \phi_t(u) W$ is a compact 
perturbation of $W$, for every $(t,u)$ in the domain of the local
flow.

\begin{lem}
\label{lemuno}
The set of smooth vector fields which satisfy (C2) with respect to
$W$ is a module over $C^{\infty}(M,\R)$.
\end{lem}

\noindent Indeed, the last term in the identity
\[
[D(\varphi X)(u),P]P = \varphi(u) [DX(u),P]P + D\varphi(u)[P\cdot]
(I-P)X(u)  
\]
has rank one.
The importance of conditions (C1) and (C2) lies in the following fact:

\begin{prop}[\cite{ama03c}, Theorem 1.6] 
\label{pop1}
Assume that the vector field $X$ has only hyperbolic rest points and
satisfies (C1), (C2). Then for every pair of rest points $x,y$, the
immersed submanifolds $W^u(x)$ and $W^s(y)$ have Fredholm intersection
of index $i(x,W) - i(y,W)$. 
\end{prop}

When the Fredholm index of these intersections is positive, the vector
field can be locally perturbed so as to make these intersections
transversal at some point. For instance, we have the following:

\begin{prop}[\cite{ama04}]
\label{fant} 
Let $x,y$ be hyperbolic rest points of a smooth vector field
$X:M\rightarrow H$ which satisfies (C1) and (C2). Let $u_0\in W^u(x;X)
\cap W^s(y;X)$ and let $u$ be the integral line of $X$ such that $u(0)=u_0$.
If $i(x,W) > i(y,W)$, for every $\epsilon>0$ and every $k\in \N$
there exists a smooth vector field $Y:M\rightarrow H$ such that:
\begin{enumerate}
\item $Y$ satisfies (C1) and (C2);
\item $Y=X$ on $M\setminus B_{\epsilon}(u_0)$ and on $u(\R)$;
\item $\|Y-X\|_{C^k} < \epsilon$;
\item $W^u(x;Y)$ and $W^s(y;Y)$ meet transversally at $u_0$.
\end{enumerate}
\end{prop}

If the vector field $X$ is the negative gradient of a
Morse function $f$ with respect to some Riemannian structure $g_0$ on
$M$, the rest points of $X$ are exactly the critical points of $f$, and
the unstable and the stable manifolds of these points are
embedded submanifolds. In this case, we can perturb the Riemannian structure
$g_0$ in such a way that the new negative gradient flow
of $f$ has the unstable and stable manifolds of pairs of critical points 
intersecting transversally:

\begin{prop}[\cite{ama04}]
\label{pop2}
Assume that $f$ is a smooth Morse function on the open subset
$M\subset H$ endowed with the Riemannian structure
$g_0$, such that the vector field $-\mathrm{grad}_{g_0}\, f$ satisfies (C1)
and (C2) with respect to a closed linear subspace $W$. Let $k\in \N$.
Then there exists a $C^k$-dense set of smooth
Riemannian structures $g$, equivalent to $g_0$,
such that $-\mathrm{grad}_g\, f$ satisfies (C1) and (C2) with respect to
$W$ and the intersections $W^u(x;f,g)\cap W^s(y;f,g)$ are
transverse, for every pair of critical points $x,y$.
\end{prop}

\section{Preliminary lemmata}

We recall that the support of a diffeomorphism $\phi$ is the closure of
the set where $\phi \neq \mathrm{id}$. 

\begin{lem}
\label{luno}
Let $A_0\in GL(H)$, $r_1>0$. Then there exists a number $r_0\in
]0,r_1[$ and a smooth diffeomorphism
$\phi: H \rightarrow H$ with support in the ball of radius $r_1$,
such that $\phi(u)=A_0 u$ for $|u|\leq r_0$. 
\end{lem}

\proof
By Kuiper theorem \cite{kui65}, the general linear group $GL(H)$ is 
connected, so we can find a
curve $A\in C^{\infty}(\R,GL(H))$ such that $A(t)=A_0$ for $t\leq
r_1/2$ and $A(t)=I$ for $t\geq r_1$. Since the continuous function
$t\mapsto \|A(t)^{-1} A^{\prime}(t)\|$ has support in $]0,r_1]$, it is
    easy to find a number $r_0\in ]0,r_1[$ and a smooth function
        $\mu:]0,+\infty[ \rightarrow ]0,+\infty[$ such that $\mu(t) =
        t$ for $t\leq r_0$ and $t\geq r_1$, and
\begin{equation}
\label{disu}
\mu^{\prime}(t) > \|A(t)^{-1} A^{\prime}(t)\| \mu(t) \quad \forall t>0.
\end{equation}
The map $\phi:H \rightarrow H$ defined by
\[
\phi(x) = \begin{cases} \frac{1}{|x|} \mu(|x|) A(|x|) x, & \mbox{if }
  x\neq 0, \\ 0, & \mbox{if } x=0, \end{cases}
\]
is clearly smooth, coincides with $x \mapsto A_0 x$ for $|x|\leq r_0$
and with the identity for $|x|\geq r_1$. 

We claim that $\phi$ is bijective. Indeed, if $v\in H$, $v\neq 0$, the
equation $\phi(x)=v$ is equivalent to
\[
x=tu, \quad \mu(t) A(t) u = v,
\]
where $t>0$ and $|u|=1$. By (\ref{disu}),
\begin{eqnarray*}
\frac{d}{dt} \left( \mu(t)^{-2} |A(t)^{-1} v|^2 \right) = -2
\mu(t)^{-2} \left( \frac{\mu^{\prime}(t)}{\mu(t)} |A(t)^{-1} v|^2 +
\langle A(t)^{-1} v,A(t)^{-1} A^{\prime}(t) A(t)^{-1} v \rangle
\right) \\ \leq -2\mu(t)^{-2} \left( \frac{\mu^{\prime}(t)}{\mu(t)} -
\|A(t)^{-1} A^{\prime}(t)\| \right) |A(t)^{-1} v|^2 < 0,
\end{eqnarray*} 
so the function $t\mapsto \mu(t)^{-1} |A(t)^{-1}v|$ is strictly
decreasing. Moreover this function maps $]0,+\infty[$ onto
$]0,+\infty[$. Therefore there exists a unique $t>0$ such that
$\mu(t)^{-1} |A(t)^{-1} v|=1$, and $x=tu=t\mu(t)^{-1} A(t)^{-1} v$ is
the unique solution of $\phi(x)=v$.

We claim that $\phi$ is a diffeomorphism. By the inverse mapping
theorem it is enough to show that $D\phi(x)$ is invertible for every
$x\in H$. This is certainly true if $x=0$, and if $x=tu$ with $t>0$
and $|u|=1$ we have
\[
D\phi(x) [v] = \frac{\mu(t)}{t} A(t) v \quad \forall v\in \langle u
\rangle^{\perp}.
\]
Therefore $D\phi(x)$ is invertible if and only if the vector
$\partial_t \phi(tu)|_{t=|x|}$ does not belong to the hyperplane
\[
D\phi(x) \langle u\rangle^{\perp} = A(t) \langle u \rangle^{\perp} =
\langle {A(t)^*}^{-1} u \rangle^{\perp}.
\]
By (\ref{disu}) the quantity
\[
\langle\partial_t \phi(tu) , {A(t)^*}^{-1} u\rangle = \langle
(\mu^{\prime} A + \mu A^{\prime}) u, {A^*}^{-1}u\rangle =
\langle (\mu^{\prime}+ \mu A^{-1}A^{\prime})u,u\rangle \geq
\mu^{\prime} - \mu \|A^{-1} A^{\prime}\|
\]
is strictly positive for every $t>0$, hence $\partial_t
\phi(tu)|_{t=|x|}$ does not belong to the hyperplane $D\phi(x) \langle
u\rangle^{\perp}$. 
\qed      

The following result is a simple addendum to the Morse Lemma:

\begin{lem}
\label{lue}
Let $f:H\rightarrow \R$ be a smooth function such that $f(0)=0$,
$Df(0)=0$, and $D^2 f(0) = A$  is invertible. For every $r_1>0$ there
exists a smooth diffeomorphism $\phi:H \rightarrow H$ with support in
the ball of radius $r_1$, such that $\phi(0)=0$, $D\phi(0)=I$, and
\begin{equation}
\label{ml}
f(u) = \frac{1}{2} \langle A \phi(u), \phi(u) \rangle,
\end{equation}
for $|u|\leq r_0$.
\end{lem}  

Indeed, Palais' proof of the Morse Lemma \cite{pal63} produces a
local diffeomorphism at $0$ $\phi$ verifying $\phi(0)=0$,
$D\phi(0) = I$, and (\ref{ml}). Let $\chi:\R \rightarrow $ be a smooth
function with compact support and such that $\chi=1$ in a neighborhood
of 0. Therefore for $\epsilon>0$ small, 
\[
\phi_{\epsilon}(u) = u + \chi\left( \frac{|u|}{\epsilon} \right)
(\phi(u) - u)  
\]
defines a smooth map on $H$, which coincides with $\phi$ in a 
neighborhood of 0, hence satisfying (\ref{ml}) therein and
$D\phi_{\epsilon}(0)=I$, and which
coincides with the identity for $|u|\geq r_1$. Since
\[
\lim_{\epsilon \rightarrow 0} \|D\phi_{\epsilon} - I\|_{\infty} = 0,
\]
$\phi_{\epsilon}$ is a global diffeomorphism for $\epsilon$ small. 

The Banach manifolds embedding theorem of Eells and Elworthy \cite{ee70}
implies that every infinite dimensional separable Hilbert manifold has
a smooth open embedding into the Hilbert space $H$. 
Here we can ask that this embedding 
satisfies also some more conditions: 

\begin{lem}
\label{redux}
Let $(M,g_0)$ be a Riemannian Hilbert manifold, let $f:M\rightarrow
\R$ be a smooth Morse function whose critical points have infinite
Morse index and co-index.
Let $V:\crit(f) \rightarrow Gr_{\infty,\infty}(H)$ be an arbitrary
function. Then there exists a smooth open embedding $\psi:M\rightarrow H$
such that for every $x\in \crit(f)$:
\begin{enumerate}
\item $D\psi(x):T_x M \rightarrow H$ is an isometry;
\item there exists $r(x)>0$ and a self-adjoint operator 
$A(x)\in GL(H)$ whose negative eigenspace is $V(x)$, and such that 
\[
f(\psi^{-1}(\psi(x) + u)) = f(x) + \frac{1}{2} \langle A(x) u,u
\rangle,
\]
for every $u\in B_{r(x)}(0)$.
\end{enumerate}
\end{lem}

\proof
Let $\psi_0:M \rightarrow H$ be an open embedding, as given by Eells
and Elworthy's theorem \cite{ee70}. For every $x\in \crit(f)$ let
$T(x)\in GL(H)$ be such that $T(x) D\psi_0(x)$ is an isometry from
$(T_x M,g_0)$ to $H$, and such that $T(x) D\psi_0(x)$ maps the negative
eigenspace of the $g_0$-Hessian of $f$ at $x$ onto $V(x)$. The
possibility of finding such a $T(x)$ follows from the fact that the
unitary group of $H$ acts transitively on the component
$Gr_{\infty,\infty}(H)$ of the Grassmannian of $H$.
By Lemma
\ref{luno}, there exists a smooth diffeomorphism $\phi_1: H
\rightarrow H$ supported in a neighborhood of $\psi_0(\crit(f))$ such
that $\phi_1(\psi_0(x)+u) = \psi_0(x) + T(x)u$ for every $x\in
\crit(f)$ and $|u|$ small. Then $f_1 = f\circ \psi_0^{-1} \circ
\phi_1^{-1}$ is a smooth Morse function on an open subset of $H$, and
the negative eigenspace of 
$A(x) = D^2 f_1(\phi_1(\psi_0(x)))$ is $V(x)$, for every $x\in
\crit(f)$. By Lemma \ref{lue}, we can find a smooth diffeomorphism
$\phi_2 : H \rightarrow H$ supported in a neighborhood of $\crit(f_1)$
such that for every $x\in \crit(f_1)$ there holds
$\phi_2(x)=x$, $D\phi_2(x)=I$, and 
\[
f_1(u) = f_1(x) + \frac{1}{2} \langle A(x) \phi_2(u), \phi_2(u)\rangle,
\]
for every $u$ close to $x$. The smooth open embedding 
$\phi_2^{-1} \circ \phi_1 \circ \psi_0$
satisfies all the requirements.
\qed  

\section{Proof of the theorem}

Fix some $W\in Gr_{\infty,\infty}(H)$.
By Lemma \ref{redux} we may assume that $M$ is an open subset of $H$,
that the Riemannian structure $g_0$ coincides with $\langle \cdot,
\cdot \rangle$ at every critical point of $f$, and that every $x\in
\crit(f)$ has a neighborhood $U_x\subset M$ such that
\begin{equation}
\label{mf}
f(u) = f(x) + \frac{1}{2} \langle A(x) (u-x),(u-x) \rangle \quad
\forall u\in U_x,
\end{equation}
where the negative eigenspace $V^-(A(x))$ of the
invertible self-adjoint operator $A(x)=D^2f(x)$ is a compact
perturbation of $W$ with $\dim (V^-(A(x)),W) = a(x)$.
Up to replacing $g_0$ by a uniformly equivalent metric and up to replacing the
neighborhoods $U_x$ by smaller ones, we may assume
that $g_0$ coincides with $\langle \cdot,\cdot \rangle$ on the whole
$U_x$, for every $x\in \crit(f)$.

We denote by $|\cdot|_{g_0(u)}$ the norm associated to the inner
product given by $g_0$ at $u$, $|v|_{g_0(u)}^2 = g_0(u)[v,v]$, and by
$\|\cdot \|_{g_0(u)}$ the corresponding operator norm on the space of
bounded linear operators on $H$,
\[
\|L\|_{g_0(u)} = \sup_{\substack{v\in H \\ |v|_{g_0(u)}=1}}
|Lv|_{g_0(u)}.
\]  

We will define a vector field $X$ on $M$, having $f$ as a Lyapunov
function, and such that the unstable and stable manifolds of its rest
points have the required properties.  
Then we will find a Riemannian structure $g_1$ for which
$X=-\mathrm{grad}_{g_1} f$. Finally we will perturb $g_1$ in order
to have a Morse-Smale flow. 

The first step of the proof consists in modifying the curves $u_i$ in order to
make them integral lines in the intersections $W^u(x_i)\cap W^s(y_i)$,
with respect to the vector field $X$ we are defining.
 
\begin{lem}
\label{step1}
There exist $\theta_0\in [0,1[$, $T>0$ and for every $i=1,\dots,n$
    there are curves $v_i\in C^{\infty}(\R,M)\cap C^0(\overline{\R},M)$ 
such that $v_i(\R) \cap v_j(\R) = \emptyset$ if $i\neq j$, 
\begin{eqnarray*}
v_i(-\infty)=x_i, \quad v_i(+\infty)=y_i, \quad v_i^{\prime}(t) = -
\mathrm{grad}_{g_0} f(v_i(t)) \mbox{ for } |t|\geq T, \\
v_i([-\infty,-T]) \subset U_{x_i}, \quad v_i([T,+\infty]) \subset U_{y_i},
\end{eqnarray*}
and
\begin{equation}
\label{edun}
|v_i^{\prime}(t) + \mathrm{grad}_{g_0} f(v_i(t))|_{g_0(v_i(t))} \leq
 \theta_0 |\mathrm{grad}_{g_0} f(v_i(t))|_{g_0(v_i(t))} \quad \forall
 t\in \R.
\end{equation}
\end{lem}

\proof
By modifying the curve $u_i:[0,1]\rightarrow M$ in a neighborhood of
$0$ and $1$, we can join it with two negative gradient flow orbits in the
unstable manifold of $x_i$ and in the stable manifold of $y_i$, 
constructing a curve $\tilde{v}_i \in
C^{\infty}(\R,M) \cap C^0(\overline{\R},M)$ satisfying all the
requirements, with the possible exception of (\ref{edun}), which will
be achieved by a reparameterization.

The construction of the curve $\tilde{v_i}$ is the following, where we
omit the index $i$ in order to simplify the notation.
Let $r>0$ be such that $B_r(x)\subset U_{x}$ and $B_r(y) \subset
U_{y}$, so that $f$ has the quadratic form (\ref{mf}) on $B_r(x)$
and on $B_r(y)$, and
\[
\mathrm{grad}_{g_0} f(\xi) = A(x) (\xi - x) \quad \forall \xi
\in B_r(x), \quad  \mathrm{grad}_{g_0} f(\xi) = A(y) (\xi - y)
\quad \forall \xi \in B_r(y).
\]
Let $\epsilon>0$ be such that $u(\epsilon)\in B_r(x)$ and
$u(1-\epsilon) \in B_r(y)$.  
Since the level sets $\Sigma^- = \set{\xi \in B_r(x)}{f(\xi) =
  f(u(\epsilon)) }$ and  $\Sigma^+ = \set{\xi \in B_r(y)}{f(\xi) =
  f(u(1-\epsilon)) }$ are connected and meet $x+V^-(A(x))$,
respectively $y+V^+(A(y))$, we can find a continuous piecewise
smooth curve $w:\R \rightarrow M$ such that $w(t) = u(t)$ for
$\epsilon \leq t \leq 1-\epsilon$, $w(t)\in \Sigma^-$ for $0\leq t
\leq \epsilon$, $w(t)\in \Sigma^+$ for $1-\epsilon \leq t
\leq 1$, $w(0)\in x + V^-(A(x))$, $w(1)\in y + V^+(A(y))$,
and
\[
w(t) = x + e^{-tA(x)} (w(0)-x) \quad \forall t<0, \quad
w(t) = y + e^{-(t-1)A(y)} (w(1)-y) \quad \forall t>1.
\]
The sets $\Sigma^-$ and $\Sigma^+$ are leafs of the codimension one
foliation given by the level sets of $f$. Therefore, it is easy to
modify $w$ in a neighborhood of $[0,\epsilon]$ and $[1-\epsilon,1]$,
obtaining a curve $\tilde{v}$ with the required properties. By
slightly perturbing these curves, we may assume that the supports of
$\tilde{v}_i$ are pairwise disjoint.

Let $\tau$ be a solution of the ODE
\[
\tau^{\prime} = - \frac{ \langle \mathrm{grad}_{g_0}
  f(\tilde{v}(\tau)) , \tilde{v}^{\prime} (\tau)
  \rangle_{g_0(\tilde{v}(\tau))}
  }{|\tilde{v}^{\prime}(\tau)|_{g_0(\tilde{v}(\tau))}^2} .
\]
Then $\tau^{\prime}(t)>0$ for every $t\in \R$, and
$\tau^{\prime}=1$ for $|t|$ large, so $\tau:\R \rightarrow \R$ is
a diffeomorphism, and it is a constant shift for $t$ in a neighborhood
of $-\infty$ and of $+\infty$. Setting $v(t) =
\tilde{v}(\tau(t))$, we have that
\begin{eqnarray*}
| v^{\prime} + \mathrm{grad}_{g_0} f(v)|_{g_0(v)} = \left| -
  \frac{ \langle \mathrm{grad}_{g_0} f(\tilde{v}(\tau)) ,
  \tilde{v}^{\prime} (\tau) \rangle_{g_0(\tilde{v}(\tau))}
  }{|\tilde{v}^{\prime}(\tau)|_{g_0(\tilde{v}(\tau))}^2}
  \, \tilde{v}^{\prime}  (\tau) +
  \mathrm{grad}_{g_0} f(\tilde{v}(\tau)) \right| \\
= \min_{\lambda\in \R} |\lambda \tilde{v}^{\prime} (\tau) +
  \mathrm{grad}_{g_0} f(\tilde{v}(\tau))
  |_{g_0(\tilde{v}(\tau))} = |\mathrm{grad}_{g_0}
  f(v)|_{g_0(v)} \sin \alpha,
\end{eqnarray*}
where $\alpha=\alpha(t)$ is the angle between the vector
$\tilde{v}^{\prime}(\tau(t))$ and $-\mathrm{grad}_{g_0}
(\tilde{v}(\tau(t)))$, with respect to the inner product
$g_0(\tilde{v}(\tau(t)))$. By the properties of $\tilde{v}$, the
continuous function $\alpha$ takes values in $[0,\pi/2[$ and it
vanishes for $|t|$ large. Hence $\alpha$ is bounded away from $\pi/2$,
and (\ref{edun}) holds with $\theta_0= \max \sin \alpha <1$.
\qed  

The second step consists in defining the vector field $X$ in a
neighborhood of the support of the curves $v_i$:

\begin{lem}
\label{step2}
Let $\theta\in ]\theta_0,1[$. There exist an
open neighborhood $U_*$ of 
$\bigcup_{i=1}^n v_i(\overline{\R})$ and a smooth
vector field $X_*:U_* \rightarrow H$ such that:
\begin{enumerate}
\item $X_*$ satisfies (C1) and (C2) with respect to $W$;
\item $X_*(u) = - \mathrm{grad}_{g_0} f(u)$ in a neighborhood of
  $\{x_1,\dots,x_n,y_1,\dots,y_n\}$;
\item $v_i^{\prime}(t) = X_*(v_i(t))$ for every $t\in \R$ and every
  $i=1,\dots,n$;
\item $|X_*(u) + \mathrm{grad}_{g_0} f(u)|_{g_0(u)} \leq \theta
  |\mathrm{grad} f(u)|_{g_0(u)}$ for every $u\in U_*$;
\item $W^u(x_i;X_*)$ and $W^s(y_i;X_*)$ intersect transversally along
  the flow line $v_i$, for every $i=1,\dots,n$.
\end{enumerate}
\end{lem}

\proof
We start by defining a smooth vector field $\tilde{X}$ on a
neighborhood $\tilde{U}$ of  $\bigcup_{i=1}^n v_i(\overline{\R})$
satisfying (i), (ii), (iii). For every $i=1,\dots, n$, let $U_i:
[-2T,2T] \rightarrow U(H) \cap (I+\mathcal{L}_c(H))$ be a smooth path
in the group of unitary operators on $H$ which are compact
perturbations of the identity, such that
\[
U_i(t) \langle v_i^{\prime}(0) \rangle^{\perp} = \langle v_i^{\prime}
(t) \rangle^{\perp} \quad \forall t\in [-2T,2T].
\]
For instance, if $P_i(t)$ denotes the orthogonal projection onto the
one-dimensional linear subspace $<v_i^{\prime}(t)>$, one can choose
$U_i$ to be the solution of the linear Cauchy problem
\[
\begin{cases}
U_i^{\prime} (t) = [P_i^{\prime}(t), P_i(t)] U_i(t), \\
U_i(0)=I. \end{cases}
\]
See \cite{kat80}, section II - \S 4.5. If $r>0$ is small enough, the map
\[
\psi_i : ]-2T,2T[ \times (\langle v_i^{\prime} (0) \rangle^{\perp}
    \cap B_r(0) ) \rightarrow H, \quad (t,\xi) \mapsto v_i(t) + U_i(t)\xi,
\]
is a diffeomorphism onto a tubular neighborhood $V_i$ of
$v_i(]-2T,2T[)$. By choosing a smaller $r$, we may also assume that
    $V_i\cap \crit(f)=\emptyset$ and that $V_i\cap V_j=\emptyset$ if
    $i\neq j$. The smooth vector field ${\psi_i}_*(\partial/\partial
    t)$,
\[
{\psi_i}_*\left( \frac{\partial}{\partial t} \right) (u) =
D\psi_i(\psi_i^{-1} (u)) [(1,0)]  \quad \forall u\in V_i, 
\]
has $v_i|_{]-2T,2T[}$ as an integral flow line and satisfies (C2) with
respect $W$, as its differential at every point is a compact
operator. Let $V_0\subset \bigcup_{i=1}^n U_{x_i} \cup U_{y_i}$ be an
open neighborhood $\bigcup_{i=1}^n v_i(\overline{\R} \setminus
  ]-2T,2T[)$ such that $v_i^{-1}(V_0) \cap [-T,T] = \emptyset$ for
    every $i=1,\dots,n$. Let $\{\chi_i\}_{i=0}^n$ be a smooth
    partition of unity subordinated to the open covering
    $\{V_i\}_{i=0}^n$ of $\tilde{U}:= V_0\cup V_1 \cup \dots \cup
    V_n$. By Lemma \ref{lemuno}, the smooth vector field 
\[
\tilde{X} : \tilde{U} \rightarrow H, \quad \tilde{X} = -\chi_0
\mathrm{grad}_{g_0} f + \sum_{i=1}^n \chi_i {\psi_i}_* \left(
\frac{\partial}{\partial t} \right),
\]  
satisfies (i), (ii), (iii). By Proposition \ref{fant}, we can modify
$\tilde{X}$ in a small neighborhood of the set $\{v_1(0),\dots,v_n(0)\}$,
obtaining a smooth vector field $X_*:\tilde{U} \rightarrow H$ which
satisfies also (v). By (ii), (iii), and (\ref{edun}), $X_*$ satisfies also
(iv) in a smaller neighborhood $U_* \subset \tilde{U}$ of
$\bigcup_{i=1}^n v_i (\overline{\R})$.
\qed

The third step consists in defining the vector field $X$ in a
neighborhood of the critical points other than $x_i$ and $y_i$.
Let $x\in \crit(f)\setminus \bigcup_{i=1}^n \{x_i,y_i\}$. 
We define $X_x$ to be the affine vector field
\[
X_x(u) = - A(x) (u-x).
\]
Its unique rest point is $x$, and the positive eigenspace of $DX_x(x)$
is $V^+(DX_x(x)) = V^+(-A(x)) =
V^-(A(x))$, which was chosen to be a compact perturbation of $W$,
so $X_x$ satisfies (C1) with respect to $W$, and the
relative Morse index of $x$ is $i(x,W) = a(x)$. Moreover by (\ref{ecom}),
\[
[DX_x(u),P]P = -[A(x),P]P = (P-I)A(x)P = (P-I)A(x) P_{V^+(A(x))} P -
(I-P) P_{V^-(A(x))} A(x) P
\]
is compact, so $X_x$ satisfies (C2) with respect to $W$. Since $g_0$
coincides with $\langle \cdot,\cdot \rangle$ on $U_x$, by (\ref{mf}) we have
\begin{equation}
\label{pg1}
X_x(u) = - \mathrm{grad}_{g_0} f(u) \quad \forall u\in U_x.
\end{equation}

The fourth step consists in defining the vector field $X$ in a
neighborhood of all the remaining points.
We shall denote by $M_0$ the complement of 
$\bigcup_{i=1}^n v_i(\overline{\R})$ in $M$.
Let $x\in M_0\setminus \crit(f) $, and consider the constant vector field
\[
X_x (u) = - \mathrm{grad}_{g_0} f(x),
\]
which trivially satisfies (C2) with respect to $W$. Let $U_x\subset
M_0\setminus \crit(f)$ be an open neighborhood of $x$ such that
\begin{equation}
\label{pg2}
|X_x (u) + \mathrm{grad}_{g_0} f(u)|_{g_0(u)} \leq \theta\,
|\mathrm{grad}_{g_0} f(u)|_{g_0(u)} \quad \forall u\in U_x.
\end{equation}

Now we can patch the vector fields $X_x$, $x\in M_0$, and $X_*$ by a
partition of unity. Indeed, let $\{\varphi_j\}_{j\in J}$ be a 
smooth partition of unity
subordinated to a locally finite refinement of the open covering
$\{U_x\}_{x\in M_0\cup \{*\}}$ of $M$: 
$\varphi_j$ is non-negative, it has support
in $U_{\omega(j)}$, $\omega(j)\in M_0 \cup \{*\}$,
every point of $M$ has a neighborhood which
intersects finitely many supports of $\varphi_j$'s, and $\sum_{j\in J}
\varphi_j =1$. Consider the smooth vector field on $M$:
\[
X(x) = \sum_{j\in J} \varphi_j(x) X_{\omega(j)} (x).
\]
By Lemma \ref{lemuno}, $X$ satisfies (C2) with respect to $W$. 
From (\ref{pg1}) and from the fact that the condition in Lemma
\ref{step2} (iv) and (\ref{pg2}) are convex conditions,
we deduce
\begin{equation}
\label{pg3}
|X (u) + \mathrm{grad}_{g_0} f(u)|_{g_0(u)} \leq 
\theta\, |\mathrm{grad}_{g_0} f(u)|_{g_0(u)} \quad
 \forall u\in M.
\end{equation}
In particular, the set of rest points of $X$ is $\crit (f)$. By Lemma
\ref{step2} (ii) and (\ref{pg1}),
\begin{equation}
\label{enc}
X(u) = - \mathrm{grad}_{g_0} f(u) \quad \forall x\in U_x^{\prime}, \;
\forall x\in \crit(f),
\end{equation}
for some neighborhood $U_x^{\prime}$ of the critical point $x$.
So $DX(x)=-A(x)$ for every critical point $x$, hence $X$ satisfies (C1) 
with respect to $W$, and the relative Morse index of $x$ is 
\begin{equation}
\label{idi}
i(x,W)= \dim (V^-(A(x)),W) = a(x).
\end{equation}
By Lemma \ref{step2} (iii) and (v), for every $i=1,\dots,n$ the
unstable manifold $W^u(x_i;X)$ meets the stable manifold $W^s(y_i;X)$
transversally along the flow line $v_i$ (indeed, by construction
$v_i(t)\in \supp \varphi_j$ if and only if $\omega(j)=*$).

We claim that there exists a Riemannian structure $g_1$ on $M$,
uniformly equivalent to $g_0$, such that $X=-\mathrm{grad}_{g_1} f$. The
construction of $g_1$ makes use of the following lemma, whose
straightforward proof is left to the reader:

\begin{lem}
Let $x\in M$ and $u,v\in H$ be such that $g_0(u,v)_x > 0$. Then the
bounded linear operator $L_x(u,v)$ defined by
\[
L_x(u,v) w = w - \frac{g_0( w,u )_x}{|u|^2_{g_0(x)}} u +
\frac{g_0(w,v)_x}{g_0(u,v)_x}v, \quad w\in H,
\]
is $g_0(x)$-self-adjoint, $g_0(x)$-positive, invertible, satisfies
\[
L_x(u,v)u=v, \quad L_x(u,u)=I,
\]
and
\begin{equation}
\label{stm} \begin{split}
\|L_x(u,v)\|_{g_0(x)} 
 \leq 1 + \frac{|v|^2_{g_0(x)}}{g_0(u,v)_x}, \\
\|L_x(u,v)^{-1}\|_{g_0(x)} 
 \leq \left( 1+ \frac{|u|_{g_0(x)} |v|_{g_0(x)}}{g_0(u,v)_x}
\right)^2 + \frac{|u|^2_{g_0(x)}}{g_0(u,v)_x}. \end{split}
\end{equation}
\end{lem}

Every Riemannian structure $g_1$ on $M$ can be written as
\[
g_1(u,v)_x = g_0(G(x)u,v)_x, \quad \forall x\in M,
\]
for some smooth map $G:M \rightarrow GL(H)$, such that $G(x)$ is
$g_0(x)$-symmetric and $g_0(x)$-positive for every $x$. 
In this case, the $g_1$-gradient of $f$ is related to the
$g_0$-gradient of $f$ by $\mathrm{grad}_{g_0} f = G\,
\mathrm{grad}_{g_1} f$. Moreover, $g_1$ is uniformly
equivalent to $g_0$ if and only if
\[
\sup_{x\in M} \|G(x)\|_{g_0(x)} < +\infty, \quad  
\sup_{x\in M} \|G(x)^{-1}\|_{g_0(x)} < +\infty.
\]
Then we can define a Riemannian structure $g_1$ by setting
\[
G(x) = \begin{cases} I & \mbox{for } x\in \bigcup_{y\in \crit(f)}
  U_y^{\prime}, \\ L_x(X(x),-\mathrm{grad}_{g_0} f(x)) & \mbox{for }
  x\in M \setminus \crit (f). \end{cases}
\]
Indeed, the above formula defines a smooth map because the map
$(x,u,v)\mapsto L_x(u,v)$ is smooth on
\[
\set{(x,u,v)\in M \times H \times H}{g_0(u,v)_x >0},
\]
and because of (\ref{enc}) and $L_x(u,u)=I$. Since $L_x(u,v)u=v$, we
obtain for $x\in M\setminus \crit(f)$,
\[
G(x)X(x) = L_x (X(x),-\mathrm{grad}_{g_0} f(x)) X(x) = -
\mathrm{grad}_{g_0} f(x) = - G(x) \mathrm{grad}_{g_1} f(x),
\]
from which $X = -\mathrm{grad}_{g_1} f$ on $M$. Finally, an easy
computation shows that (\ref{pg3}) and (\ref{stm}) imply
\[
\|G(x)\|_{g_0(x)} \leq 1 + \frac{1}{1-\theta}, \quad
\|G(x)^{-1}\|_{g_0(x)} \leq \frac{4}{(1-\theta)^2} +
\frac{(1+\theta)^2}{1-\theta}, \quad \forall x\in M,
\]
so $g_1$ is uniformly equivalent to $g_0$.

By Proposition \ref{pop2}, we can perturb $g_1$ obtaining a Riemannian 
structure $g$ on
$M$ such that the vector field $-\mathrm{grad}_g f$ still satisfies (C1) and
(C2) with respect to $W$, all the intersections between unstable
and stable manifolds are transverse, and the transverse intersections
at $v_i$ are preserved, for every $i=1,\dots,n$. 
Then Proposition \ref{pop1} and (\ref{idi}) imply that
\[
\dim W^u(x;f,g) \cap W^s(y;f,g) = i(x,W) - i(y,W) = a(x) - a(y),    
\]
whenever such intersection is non-empty.

\providecommand{\bysame}{\leavevmode\hbox to3em{\hrulefill}\thinspace}
\providecommand{\MR}{\relax\ifhmode\unskip\space\fi MR }
\providecommand{\MRhref}[2]{%
  \href{http://www.ams.org/mathscinet-getitem?mr=#1}{#2}
}
\providecommand{\href}[2]{#2}

\end{document}